\documentclass[11pt, reqno]{amsart}
\usepackage[T1]{fontenc}
\usepackage{amsfonts}
\usepackage{amsmath}
\usepackage{amsthm}
\usepackage{amssymb}
\usepackage[textwidth=15.5cm, textheight=20.5cm, hmarginratio=1:1, vmarginratio=1:1]{geometry}
\usepackage{graphicx}
\usepackage{xcolor}
\usepackage{mathtools}
\usepackage[colorlinks=true, linkcolor=blue]{hyperref}
\usepackage{cleveref}
\usepackage[english]{babel}
\usepackage{lineno}
\usepackage{float}

\newtheorem{theorem}{Theorem}[section]

\newtheorem{corollary}[theorem]{Corollary}
\newtheorem{lemma}[theorem]{Lemma}
\newtheorem{conjecture}[theorem]{Conjecture}
\newtheorem{proposition}[theorem]{Proposition}

\usepackage{tikz}
\usetikzlibrary{positioning}
\usetikzlibrary{decorations,arrows}
\usetikzlibrary{decorations.markings}
\numberwithin{equation}{section}

\newcommand{\Z}{\mathbb Z}

\usepackage{rustic}
\usepackage[T1]{fontenc}

\emergencystretch=\maxdimen
\title[Edge-coloring spectra of odd-fold graph covers]{Edge-coloring spectra of odd-fold graph covers and counterexamples to a conjecture}

\author{Runze Wang}
\address[]{Department of Mathematics and Computer Science, Augustana College, Rock Island, IL 61201, USA}
\email{runze.w@hotmail.com}
\thanks{}

\begin{document}

\sloppy

\begin{abstract}
For a loopless multigraph $G$ and an odd integer $\ell\ge3$, let $\Delta(G)$ be the maximum degree of $G$, let $s_\ell(G)$ be the minimum chromatic index of a cyclic $\ell$-fold cover of $G$, and let 
\begin{align*}
L_\ell(G):=\max\left\{
\Delta(G),\ 
\max_{\substack{S\subseteq V(G)\\ |S|\ge3\text{ odd}}}
\left\lceil\frac{2\ell|E(G[S])|}{\ell|S|-1}\right\rceil
\right\}.
\end{align*}
In this paper, we prove the sharp upper bound
\begin{align*}
s_\ell(G)\le\max\{\Delta(G)+1,L_\ell(G)\}.
\end{align*}
For every nonempty bridgeless cubic graph $H$, we construct two $3\ell$-regular graphs, one a multigraph and the other a simple graph. In each case, the minimum chromatic indices over cyclic and over all $\ell$-fold covers both equal the chromatic index of the multigraph obtained by replacing each edge of $H$ with $\ell$ parallel edges. Taking $H$ to be the Petersen graph gives both multigraph and simple graph counterexamples to a conjecture of Chen, van der Holst, Luo, and Ma, which asserts that $s_\ell(G)=L_\ell(G)$ always holds. When $G$ is connected and $L_\ell(G)>\Delta(G)$, we determine the common spectrum of its connected cyclic and arbitrary $\ell$-fold covers. We give polynomial-time algorithms to compute $L_\ell(G)$ and the spectrum endpoints. For fixed $\ell$, we also construct in polynomial time a connected cyclic cover with any prescribed chromatic index in the spectrum.
\end{abstract}

\vspace*{-1em}
\maketitle

\section{Introduction}

All graphs are finite and loopless, and parallel edges are permitted. Two parallel edges form a cycle of length $2$. For a graph $G=(V(G),E(G))$, a proper $k$-edge-coloring of $G$ assigns one of $k$ colors to each edge so that incident edges receive distinct colors. The chromatic index of $G$, denoted by $\chi'(G)$, is the smallest $k$ such that $G$ admits a proper $k$-edge-coloring. For a vertex set $U\subseteq V(G)$, denote by $G[U]$ the subgraph of $G$ induced by $U$. Denote by $\Delta(G)$ the maximum degree of $G$. Since each color class is a matching, we have
\begin{align*}
\chi'(G)\ge\max\{\Delta(G),\lceil\Gamma(G)\rceil\},
\end{align*}
where
\begin{align*}
\Gamma(G):=\max_{\substack{S\subseteq V(G)\\ |S|\ge3\text{ odd}}}
\frac{2|E(G[S])|}{|S|-1}.
\end{align*}
We call $\frac{2|E(G[S])|}{|S|-1}$ the density of $S$. We take $\Gamma(G)=0$ if there is no odd set of size at least $3$. The Goldberg-Seymour Theorem, proven by Chen, Jing, and Zang~\cite{chen2025goldberg} (see also~\cite{chen2024short}), states that
\begin{align}\label{eq:GS}
\chi'(G)\le\max\{\Delta(G)+1,\lceil\Gamma(G)\rceil\}.
\end{align}

For a positive integer $\ell$, an $\ell$-fold cover of $G$ is a graph $J$ with an incidence-preserving projection $\pi:J\longrightarrow G$ such that each vertex $v\in V(G)$ has a fiber $\pi^{-1}(v)$ of size $\ell$, and each base edge $uv$ lifts to a perfect matching between $\pi^{-1}(u)$ and $\pi^{-1}(v)$. We call $\pi$ the covering projection. For $S\subseteq V(G)$, its full preimage is the vertex set $\pi^{-1}(S)$. Labeling each fiber by the same $\ell$ elements identifies the matching over each oriented base edge with a permutation of the labels. Gross~\cite{gross1974voltage} introduced voltage graphs, and Gross and Tucker~\cite{gross1977coverings} showed that every graph cover can be obtained from permutation voltages.

In a cyclic $\ell$-fold cover, the permutations are translations of the cyclic group $\Z_\ell$. Fix a reference orientation of $G$ and let $\alpha:E(G)\longrightarrow\Z_\ell$ be a voltage assignment. The resulting cover $\widetilde G(\alpha,\ell)$ has vertex set $V(G)\times\Z_\ell$. Each edge $e=uv$ oriented from $u$ to $v$ lifts to the edges $(u,i)(v,i+\alpha(e))$ for $i\in\Z_\ell$. The covering projection is given on vertices by $\pi(v,i)=v$. Each lifted edge projects to the corresponding base edge $e$. Thus the full preimage of $S\subseteq V(G)$ is $S\times\Z_\ell$. The second coordinate $i$ is the sheet label, and the vertices with a fixed sheet label form a sheet. Reversing the orientation of a base edge and negating its voltage gives the same cover. In any $\ell$-fold cover $J$ of $G$, each vertex has the degree of its image, and a proper edge-coloring of $G$ lifts to a proper edge-coloring of $J$. Hence
\begin{align*}
\Delta(G)=\Delta(J)\le\chi'(J)\le\chi'(G).
\end{align*}

Edge-colorings of graph coverings have been studied in several settings. Pisanski, Shawe-Taylor, and Vrabec~\cite{pisanski1983bundles} gave sufficient conditions for a graph bundle to have chromatic index equal to its maximum degree. Plachta~\cite{plachta2021coverings} studied conditions under which a finite cover of a $3$-edge-uncolorable cubic graph remains $3$-edge-uncolorable. More recently, Luo, Ma, Stephens, and Ye~\cite{luo2026flows} proved that every bridgeless cubic graph has a $3$-edge-colorable cyclic $\ell$-fold cover for every odd $\ell\ge3$.

For a family $\mathcal{F}$ of covers, the \emph{edge-coloring spectrum} of $\mathcal{F}$, or simply the \emph{spectrum} of $\mathcal{F}$, is the set $\{\chi'(J):J\in\mathcal{F}\}$. Chen, van der Holst, Luo, and Ma~\cite{chen2026spectra} studied the spectrum of all cyclic $\ell$-fold covers of a fixed graph. They proved that, when $\ell$ is even, this spectrum contains every integer between $\Delta(G)$ and $\chi'(G)$. For odd $\ell$, they determined the spectrum for edge-critical graphs $G$ with $\Delta(G)\ge3$ and $\chi'(G)\ge\Delta(G)+1$. A graph is \emph{edge-critical} if every proper subgraph has smaller chromatic index. Their density bound~\cite[Theorem~2.1]{chen2026spectra} implies that every cyclic cover $\widetilde G$ of $G$ satisfies
\begin{align*}
\lceil\Gamma(\widetilde G)\rceil\le\max\{\Delta(G),\lceil\Gamma(G)\rceil\}.
\end{align*}
For an odd cycle $G$, every odd-fold cover has chromatic index $3$, since it is $2$-regular and has odd order.

For an integer $\ell\ge2$, let $s_\ell(G)$ be the minimum chromatic index over cyclic $\ell$-fold covers of $G$, and let $a_\ell(G)$ be the minimum chromatic index over all $\ell$-fold covers of $G$. For a positive integer $\ell$, define
\begin{align*}
L_\ell(G):=\max\left\{
\Delta(G),\ 
\max_{\substack{S\subseteq V(G)\\ |S|\ge3\text{ odd}}}
\left\lceil\frac{2\ell|E(G[S])|}{\ell|S|-1}\right\rceil
\right\},
\end{align*}
where the inner maximum is $0$ if there is no odd set of size at least $3$. Let $J$ be an $\ell$-fold cover of $G$ with covering projection $\pi$. If $\ell\ge3$ is odd and $S\subseteq V(G)$ has odd size, then $\pi^{-1}(S)$ has $\ell|S|$ vertices and induces $\ell|E(G[S])|$ edges. Each color class contains at most $\frac{\ell|S|-1}{2}$ of these edges. Together with the degree bound, this gives $a_\ell(G)\ge L_\ell(G)$.

Chen, van der Holst, Luo, and Ma~\cite[Conjecture~1]{chen2026spectra} proposed the following conjecture.

\begin{conjecture}\label{conj:chen}
Let $G$ be a loopless graph and let $\ell\ge3$ be odd. Then $s_\ell(G)=L_\ell(G)$.
\end{conjecture}

We disprove Conjecture~\ref{conj:chen} by constructing both multigraph and simple graph counterexamples for every odd $\ell\ge3$. We first prove the following sharp bound.

\begin{theorem}\label{thm:bound}
For every graph $G$ and every odd integer $\ell\ge3$,
\begin{align*}
L_\ell(G)\le a_\ell(G)\le s_\ell(G)\le\max\{\Delta(G)+1,L_\ell(G)\}.
\end{align*}
\end{theorem}

By Theorem~\ref{thm:bound}, we have $a_\ell(G)=s_\ell(G)=L_\ell(G)$ whenever $L_\ell(G)>\Delta(G)$. When $L_\ell(G)=\Delta(G)$, both $a_\ell(G)$ and $s_\ell(G)$ can exceed $L_\ell(G)$ by $1$, even for simple $3$-edge-connected graphs. We prove this by relating the chromatic indices of the covers to perfect matchings of a bridgeless cubic graph.

For a positive integer $k$, let $kH$ denote the multigraph obtained by replacing every edge of $H$ with $k$ parallel edges. If $H$ is cubic, a proper $3k$-edge-coloring of $kH$ is equivalent to a collection of $3k$ perfect matchings of $H$, with repetitions allowed, in which every edge occurs exactly $k$ times. The next theorem realizes $\chi'(\ell H)$ as the minimum chromatic index of covers of two graphs constructed from $H$.

\begin{theorem}\label{thm:main}
Let $H$ be a nonempty bridgeless cubic graph, let $\ell\ge3$ be odd, and let $q=3\ell$. There exist a $q$-regular multigraph $F_\ell(H)$ of order $3|V(H)|$ and a $q$-regular simple graph $G_\ell(H)$ of order $(q+2)|V(H)|$ such that, for each $X\in\{F_\ell(H),G_\ell(H)\}$,
\begin{align*}
L_\ell(X)&=q,\\
a_\ell(X)=s_\ell(X)&=\chi'(\ell H)\in\{q,q+1\}.
\end{align*}
If $H$ is connected, then both graphs are connected, and if $H$ is $3$-edge-connected, then both graphs are $3$-edge-connected.
\end{theorem}

By Theorem~\ref{thm:main}, either constructed graph has a $q$-edge-colorable $\ell$-fold cover if and only if $H$ has a collection of $q$ perfect matchings in which every edge occurs exactly $\ell$ times.

\begin{corollary}\label{cor:counterexamples}
For every odd integer $\ell\ge3$, there exist two $3$-edge-connected $3\ell$-regular graphs, a multigraph $F_\ell$ on $30$ vertices and a simple graph $G_\ell$ on $10(3\ell+2)$ vertices, such that, for each $X\in\{F_\ell,G_\ell\}$,
\begin{align*}
L_\ell(X)=3\ell
\quad\text{and}\quad
a_\ell(X)=s_\ell(X)=3\ell+1.
\end{align*}
Every $\ell$-fold cover of $F_\ell$ has chromatic index at least $3\ell+1$, and every $\ell$-fold cover of $G_\ell$ has chromatic index exactly $3\ell+1$.
\end{corollary}

Taking $H$ to be the Petersen graph in Theorem~\ref{thm:main} gives the graphs in Corollary~\ref{cor:counterexamples}. These graphs disprove Conjecture~\ref{conj:chen} and show that Theorem~\ref{thm:bound} is sharp for every odd $\ell\ge3$.

Let $G$ be a connected graph containing a cycle. For $\ell\ge2$, let $\mathcal{S}_\ell^{\mathrm{cyc}}(G)$ be the spectrum of connected cyclic $\ell$-fold covers of $G$, namely the set of their chromatic indices. Let $\mathcal{S}_\ell^{\mathrm{all}}(G)$ be the corresponding spectrum of all connected $\ell$-fold covers of $G$. For integers $p\le q$, let $[p,q]=\{r\in\Z:p\le r\le q\}$.

\begin{theorem}\label{thm:connected}
Let $G$ be a connected graph containing a cycle, let $\ell\ge2$, and let $k=\chi'(G)$.
\begin{enumerate}
\item[\textup{(i)}] The spectra $\mathcal{S}_\ell^{\mathrm{cyc}}(G)$ and $\mathcal{S}_\ell^{\mathrm{all}}(G)$ are intervals of integers with lower endpoints $s_\ell(G)$ and $a_\ell(G)$, respectively.
\item[\textup{(ii)}] If $k\ge\Delta(G)+2$, both spectra have upper endpoint
\begin{align}\label{eq:upper-end}
M=
\begin{cases}
k, & \text{if $\chi'(G-e)=k$ for some nonbridge $e$},\\
k-1, & \text{if $\chi'(G-e)=k-1$ for every nonbridge $e$}.
\end{cases}
\end{align}
\item[\textup{(iii)}] Assume that $\ell\ge3$ is odd. If $L_\ell(G)>\Delta(G)$, or if $k\ge\Delta(G)+2$ and $\chi'(G-e)=k-1$ for every nonbridge $e$, then both spectra are $[L_\ell(G),M]$. Here $M$ is given by \eqref{eq:upper-end} for $k\ge\Delta(G)+2$ and $M=k$ otherwise.
\end{enumerate}
\end{theorem}

In Section~\ref{sec:bound}, we prove Theorem~\ref{thm:bound} using a cyclic cover obtained by pairing parallel edges. Section~\ref{sec:constructions} gives the two graph constructions and proves Theorem~\ref{thm:main}. In Section~\ref{sec:petersen}, we prove Corollary~\ref{cor:counterexamples} by counting spokes in perfect matchings of the Petersen graph. In Section~\ref{sec:connected}, we prove Theorem~\ref{thm:connected}, establishing the interval property by changing the matching over one base edge at a time while preserving connectivity. Section~\ref{sec:algorithm} gives polynomial-time algorithms to compute $L_\ell(G)$ and the endpoints in parts~(ii) and~(iii) of Theorem~\ref{thm:connected}. For fixed odd $\ell\ge3$ with $L_\ell(G)>\Delta(G)$, we also construct in polynomial time a connected cyclic $\ell$-fold cover with any prescribed chromatic index in the common spectrum.



\section{Proof of Theorem~\ref{thm:bound}}\label{sec:bound}

For a graph $J$ and $W\subseteq V(J)$, let $\partial_J(W)$ be the set of edges with exactly one endpoint in $W$. The sets $W$ and $V(J)\setminus W$ are the two shores of this cut.

\begin{lemma}\label{lem:pairs}
For every $\ell\ge2$ and every graph $G$, there is a cyclic $\ell$-fold cover $J$ such that
\begin{align*}
2|E(J[W])|\le (\Delta(G)+1)(|W|-1)
\end{align*}
for every $W\subseteq V(J)$ that meets a fiber without containing it.
\end{lemma}

\begin{proof}
For each pair of vertices of $G$, partition the edges joining them into pairs, leaving at most one edge unpaired. Orient the edges in each pair in the same direction and assign them voltages $0$ and $1$. Following a lift of the voltage-$1$ edge and then a reverse lift of the voltage-$0$ edge changes the sheet label by $1$. Hence the lifts form a single cycle of length $2\ell$ on the two fibers. Give each unpaired edge voltage $0$. The unpaired edges form a simple spanning subgraph $R$ of $G$, and their lifts form a simple spanning subgraph $\widetilde R$ of $J$.

Let $u$ be a vertex of $G$ whose fiber meets both $W$ and its complement, and choose a vertex $x\in W$ over $u$. Each of the chosen pairs of parallel edges incident with $u$ lifts to a cycle containing the fiber over $u$ and thus contributes at least two edges to $\partial_J(W)$. The cycles are edge-disjoint, so together they contribute at least $\deg_G(u)-\deg_R(u)$ edges to the cut. Since $\widetilde R$ is simple, at most $|W|-1$ of its edges incident with $x$ have their other endpoint in $W$. The remaining edges of $\widetilde R$ incident with $x$ contribute at least $\deg_R(u)-(|W|-1)$ edges to the cut. The contributions from paired and unpaired edges are disjoint, so
\begin{align*}
|\partial_J(W)|\ge \deg_G(u)-\deg_R(u)
+\deg_R(u)-(|W|-1)=\deg_G(u)-|W|+1.
\end{align*}
Using the degree of $x$ and bounding all other degrees by $\Delta(G)$ gives
\begin{align*}
2|E(J[W])|
&=\sum_{w\in W}\deg_J(w)-|\partial_J(W)|\\
&\le\deg_G(u)+(|W|-1)\Delta(G)-\bigl(\deg_G(u)-|W|+1\bigr)\\
&=(\Delta(G)+1)(|W|-1).
\end{align*}
\end{proof}

\begin{proof}[Proof of Theorem~\ref{thm:bound}]
As shown in the introduction, $L_\ell(G)\le a_\ell(G)$. The inequality $a_\ell(G)\le s_\ell(G)$ follows from the definitions. For the upper bound $s_\ell(G)\le\max\{\Delta(G)+1,L_\ell(G)\}$, let $J$ be the cover from Lemma~\ref{lem:pairs}. Let $W\subseteq V(J)$ be odd with $|W|\ge3$. If $W$ meets a fiber without containing it, Lemma~\ref{lem:pairs} bounds its density by $\Delta(G)+1$. Otherwise $W$ is the full preimage of a set $S\subseteq V(G)$. Since $|W|=\ell|S|$ and $\ell$ is odd, $S$ is odd. Its preimage has density
\begin{align*}
\frac{2|E(J[W])|}{|W|-1}=\frac{2\ell|E(G[S])|}{\ell|S|-1}\le L_\ell(G).
\end{align*}
For a singleton $S$, the ratio is $0$ because $G$ has no loops. It follows from \eqref{eq:GS} that $\chi'(J)\le\max\{\Delta(G)+1,L_\ell(G)\}$. As mentioned in the introduction, this bound is sharp because of the graphs in Corollary~\ref{cor:counterexamples}.
\end{proof}

\section{Proof of Theorem~\ref{thm:main}}\label{sec:constructions}

\subsection{The base graphs}\label{subsec:bases}
Let $H$ be a nonempty bridgeless cubic graph and let $\ell\ge3$ be odd. Let $q=3\ell$ and $p=\frac{q-1}{2}$. In both constructions, we replace each vertex $v$ of $H$ with a graph $Q_v$ of odd order. We refer to each such replacement graph as a block. Each block has three attachment vertices $x_{v,e}$, one for each edge $e$ incident with $v$. Each attachment vertex has degree $q-1$ in the block, and every other vertex has degree $q$. For each edge $e=uv$ of $H$, add the connecting edge $f_e=x_{u,e}x_{v,e}$. The connecting edges form a matching, and the resulting graph is $q$-regular.

\textbf{The multigraph $F_\ell(H)$:}
For each $v\in V(H)$, the block $Q_v$ is a triangle with $p$ parallel edges on each side. Its three vertices are the attachment vertices. Hence $F_\ell(H)$ has $3|V(H)|$ vertices.

\textbf{The simple graph $G_\ell(H)$:}
For each $v\in V(H)$, let $A_v$ be the set of three attachment vertices, and let $B_v$ be a set of $q-1$ additional vertices. Join every vertex of $A_v$ to every vertex of $B_v$. The subgraph induced by $B_v$ is a complete graph with one perfect matching removed. The vertices of $A_v$ are pairwise nonadjacent and have degree $q-1$ in $Q_v$. Each vertex of $B_v$ has $q-3$ neighbors in $B_v$ and three neighbors in $A_v$. Each block has $q+2$ vertices, and the connecting edges form a matching between distinct blocks. Thus $G_\ell(H)$ is simple and has $(q+2)|V(H)|$ vertices.

The simple block has the smallest possible odd order for these degree requirements. An attachment vertex of degree $q-1$ requires at least $q$ vertices in the block. A block on $q$ vertices would have $q-3$ vertices of degree $q$, which is impossible in a simple graph of order $q$. As $q$ is odd, the order is therefore at least $q+2$.

Let $X$ be either $F_\ell(H)$ or $G_\ell(H)$. If $S\subseteq V(X)$ meets a block $Q_v$ without containing it, then $\partial_X(S)$ contains at least three edges of $Q_v$. For a multigraph block, the cut contains $2p=q-1$ edges of the block. For a simple block, if $A_v$ meets both sides of the cut, each vertex in $B_v$ is incident with an edge of $\partial_X(S)$ joining it to $A_v$. If $A_v$ lies on one side, some vertex in $B_v$ lies on the other and its three edges to $A_v$ cross the cut. If $S$ is a union of whole blocks and $U=\{v\in V(H): V(Q_v)\subseteq S\}$, then
\begin{align*}
|\partial_X(S)|=|\partial_H(U)|.
\end{align*}
These observations prove the connectivity assertions of Theorem~\ref{thm:main}.

If $|S|$ is odd and $S$ is a union of whole blocks, then $|U|$ is odd because every block has odd order. In a cubic graph, a cut has the same parity as the size of either shore. Since $H$ is bridgeless, a cut with an odd shore has size at least $3$. Thus $|\partial_X(S)|\ge3$ for every odd set $S$. It follows that $L_\ell(X)=q$, since $\Delta(X)=q$ and
\begin{align*}
2\ell|E(X[S])|
&=\ell\bigl(q|S|-|\partial_X(S)|\bigr)\\
&\le q(\ell|S|-1).
\end{align*}

\subsection{The lower bound for every cover}

\begin{lemma}\label{lem:lower}
The chromatic index of $\ell H$ is either $q$ or $q+1$. Every $\ell$-fold cover of $F_\ell(H)$ or $G_\ell(H)$ has chromatic index at least $\chi'(\ell H)$.
\end{lemma}

\begin{proof}
For each odd set $U\subseteq V(H)$ with $|U|\ge3$, the inequality $|\partial_H(U)|\ge3$ gives
\begin{align*}
2|E((\ell H)[U])|
=\ell\bigl(3|U|-|\partial_H(U)|\bigr)
\le q(|U|-1).
\end{align*}
Hence $\Gamma(\ell H)\le q$. Since $\Delta(\ell H)=q$, the Goldberg-Seymour Theorem gives $\chi'(\ell H)\in\{q,q+1\}$.

Let $X\in\{F_\ell(H),G_\ell(H)\}$, and assume that an $\ell$-fold cover $J$ of $X$ with projection $\pi$ admits a proper $q$-edge-coloring. As $J$ is $q$-regular, each color class is a perfect matching. For each $v\in V(H)$, let $D_v=\pi^{-1}(V(Q_v))$. The set $D_v$ has odd order, and its cut consists of the $\ell$ lifts of each of the three connecting edges at $Q_v$. Thus $|\partial_J(D_v)|=q$. Each color class meets $\partial_J(D_v)$ in an odd number of edges, since it covers the vertices of $D_v$ in pairs except at edges crossing the cut. There are $q$ colors, so each occurs exactly once on the cut.

For each color $c$, let $M_c$ consist of the edges $e$ of $H$ for which some lift of $f_e$ has color $c$. The cut condition at each $D_v$ implies that exactly one edge of $M_c$ is incident with $v$. Thus $M_c$ is a perfect matching of $H$. The $\ell$ lifts of each connecting edge $f_e$ have distinct colors because they all lie in the same cut $\partial_J(D_v)$ for an endpoint $v$ of $e$. Hence $e$ belongs to exactly $\ell$ of the matchings $M_c$, and these matchings give a proper $q$-edge-coloring of $\ell H$. Since every cover of $X$ has maximum degree $q$, this proves the claimed lower bound.
\end{proof}

\subsection{Cyclic covers attaining the lower bound}

In an edge-coloring with colors from $\mathcal{C}$, a color $c\in\mathcal{C}$ is missing at a vertex $x$ if no edge incident with $x$ has color $c$.

\begin{lemma}\label{lem:blocks}
Each type of block has a cyclic $\ell$-fold cover with the following property. For any set $\mathcal{C}$ of $q$ colors and any bijection from the lifted attachment vertices to $\mathcal{C}$, the cover has a proper edge-coloring with colors from $\mathcal{C}$ in which each lifted attachment vertex misses exactly its assigned color.
\end{lemma}

\begin{proof}
For either block it suffices to give one proper $q$-edge-coloring in which the colors missing at the lifted attachment vertices are pairwise distinct. Any specified bijection can then be obtained by permuting the colors.

\textbf{The multigraph block:}
Orient the sides of the triangle cyclically. Give the $p$ parallel edges on one side voltage $1$ and all other edges voltage $0$. Traversing the triangle once increases the sheet label by $1$. Thus the lifted block is obtained from a cycle of length $q=3\ell$ by replacing each edge with $p$ parallel edges.

For each vertex $z$ of the underlying cycle, deleting $z$ leaves a path on $q-1$ vertices. Let $N_z$ be its perfect matching. By rotational symmetry every cycle edge occurs in the same number of matchings $N_z$. This number is $p$, since there are $q$ cycle edges and $q$ matchings, each of size $p$. Assign the $p$ copies of each cycle edge bijectively to the matchings that contain it, then give the resulting $q$ matchings distinct colors. The color assigned to $N_z$ is the only color missing at $z$.

\textbf{The simple block:}
Let $A$ be the set of attachment vertices in the block, and let $B$ be the set of other vertices in the block. Let $B_0,B_1,\dots,B_{\ell-1}$ partition $B$ so that $|B_j|=3$ for $0\le j\le\ell-2$ and $|B_{\ell-1}|=2$. Give every edge with both endpoints in $B$ voltage $0$. For each $b\in B_j$ and $a\in A$, orient $ba$ from $b$ to $a$ and give it voltage $j$. For each $i\in\Z_\ell$, let $A^i=A\times\{i\}$, $B^i=B\times\{i\}$, and $B_j^i=B_j\times\{i\}$. Let $\{\mathcal{C}_i:i\in\Z_\ell\}$ be a partition of $\mathcal{C}$ into sets of size $3$. Sheet labels and the indices of the sets $B_j$ and $\mathcal{C}_i$ are read modulo $\ell$.

For fixed $i$ and $j$, the lifted edges between $B_j^i$ and $A^{i+j}$ form a complete bipartite graph with parts of sizes $|B_j|$ and $3$. By K\H{o}nig's line coloring theorem, this graph has a proper edge-coloring with the three colors in $\mathcal{C}_i$. For fixed $i$, the subgraphs indexed by $j$ are vertex-disjoint, and every vertex of $B^i$ is incident with edges of all three colors in $\mathcal{C}_i$.

Fix $h\in\Z_\ell$. The vertices in $A^h$ receive edges from $B_{h-i}^i$ for each $i$, with distinct color sets $\mathcal{C}_i$. All three colors in $\mathcal{C}_i$ occur at each vertex of $A^h$ unless $h=i+\ell-1$, when the corresponding graph is $K_{2,3}$. In that graph, each vertex in $B_{\ell-1}^i$ is incident with edges of all three colors, so each color is used at two distinct vertices of $A^h$ and is missing at the third. For each $h$ there is exactly one such $i$, and for each $i$ exactly one such $h$. Thus each lifted attachment vertex misses exactly one color and each color is missing at exactly one lifted attachment vertex.

It remains to color the edges within each $B^i$. Since $|B|$ is even, the complete graph on $B$ has a decomposition into perfect matchings. By relabeling its vertices, we may choose a decomposition containing the perfect matching deleted from the complete graph on $B$. The remaining $q-3$ matchings in $B^i$ can therefore be colored with $\mathcal{C}\setminus\mathcal{C}_i$. This completes the coloring of the lifted block.
\end{proof}

\begin{proof}[Proof of Theorem~\ref{thm:main}]
Subsection~\ref{subsec:bases} establishes the assertions about the two base graphs and $L_\ell$. By Lemma~\ref{lem:lower} it remains to find a cyclic $\ell$-fold cover of each base graph with chromatic index at most $\chi'(\ell H)$.

Let $X\in\{F_\ell(H),G_\ell(H)\}$. On each block of $X$, use the voltage assignment from Lemma~\ref{lem:blocks}, and give all connecting edges voltage $0$. Let $\widetilde X$ be the resulting cover. We may color all lifted blocks with the same set of $q$ colors and use one additional color on the lifted connecting edges, which form a matching. Thus $\chi'(\widetilde X)\le q+1$.

Now assume that $\ell H$ has a proper edge-coloring with a set $\mathcal{C}$ of $q$ colors. For each $e\in E(H)$, let $\{a_{e,i}:i\in\Z_\ell\}$ be the set of colors on the $\ell$ copies of $e$, indexed in any order. At each vertex $v$ of $H$, the three sets corresponding to the edges incident with $v$ partition $\mathcal{C}$. Thus, in the lifted block at $v$, prescribing $a_{e,i}$ as the missing color at $(x_{v,e},i)$ defines a bijection from its lifted attachment vertices to $\mathcal{C}$. Apply Lemma~\ref{lem:blocks} in every lifted block with these prescriptions. For each $e=uv$, assign $a_{e,i}$ to the lifted edge $(x_{u,e},i)(x_{v,e},i)$. This color is missing at both endpoints, so the resulting coloring of $\widetilde X$ is proper.

Hence $\chi'(\widetilde X)\le\chi'(\ell H)$ in both cases. Together with Lemma~\ref{lem:lower}, this gives $a_\ell(X)=s_\ell(X)=\chi'(\ell H)$. The voltage assignment for $\widetilde X$ depends only on the block construction.
\end{proof}

\section{Petersen counterexamples}\label{sec:petersen}
Let $P$ be the Petersen graph with vertices $a_i,b_i$ and edges
\begin{align*}
a_i a_{i+1},\quad b_i b_{i+2},\quad a_i b_i
\end{align*}
for $i\in\Z_5$. We call the vertices $a_i$ the outer vertices and the vertices $b_i$ the inner vertices. The edges $a_i b_i$ are called the spokes. The graph $P$ is $3$-edge-connected. Indeed, a cut that splits both $5$-cycles contains at least four cycle edges. A cut that splits exactly one of them contains at least two cycle edges and at least one spoke. The cut separating the two cycles consists of all five spokes.

\begin{proof}[Proof of Corollary~\ref{cor:counterexamples}]
Let $\ell\ge3$ be odd. Suppose that $\ell P$ has a proper $3\ell$-edge-coloring. Its color classes give $3\ell$ perfect matchings of $P$ in which every edge occurs exactly $\ell$ times.

Every perfect matching of $P$ contains an odd number of spokes, since the outer vertices that are incident with none of the selected spokes must be matched in pairs. A perfect matching cannot contain exactly three spokes. The indices of the two remaining outer vertices must differ by $\pm1$ modulo $5$, while the same indices on the inner cycle must differ by $\pm2$. Thus each perfect matching contains either one or five spokes.

Let $x$ of the $3\ell$ perfect matchings contain all five spokes. Counting occurrences of spokes gives
\begin{align*}
5\ell=5x+(3\ell-x)=3\ell+4x.
\end{align*}
It follows that $x=\frac{\ell}{2}$, contrary to the assumption that $\ell$ is odd. Lemma~\ref{lem:lower} now gives $\chi'(\ell P)=3\ell+1$. Apply Theorem~\ref{thm:main} with $F_\ell=F_\ell(P)$ and $G_\ell=G_\ell(P)$. A cover of a simple graph is simple. Since every cover of $G_\ell$ is also $3\ell$-regular, Vizing's theorem gives the upper bound $3\ell+1$.
\end{proof}

\begin{figure}[H]
    \centering
    \includegraphics[width=0.95\linewidth]{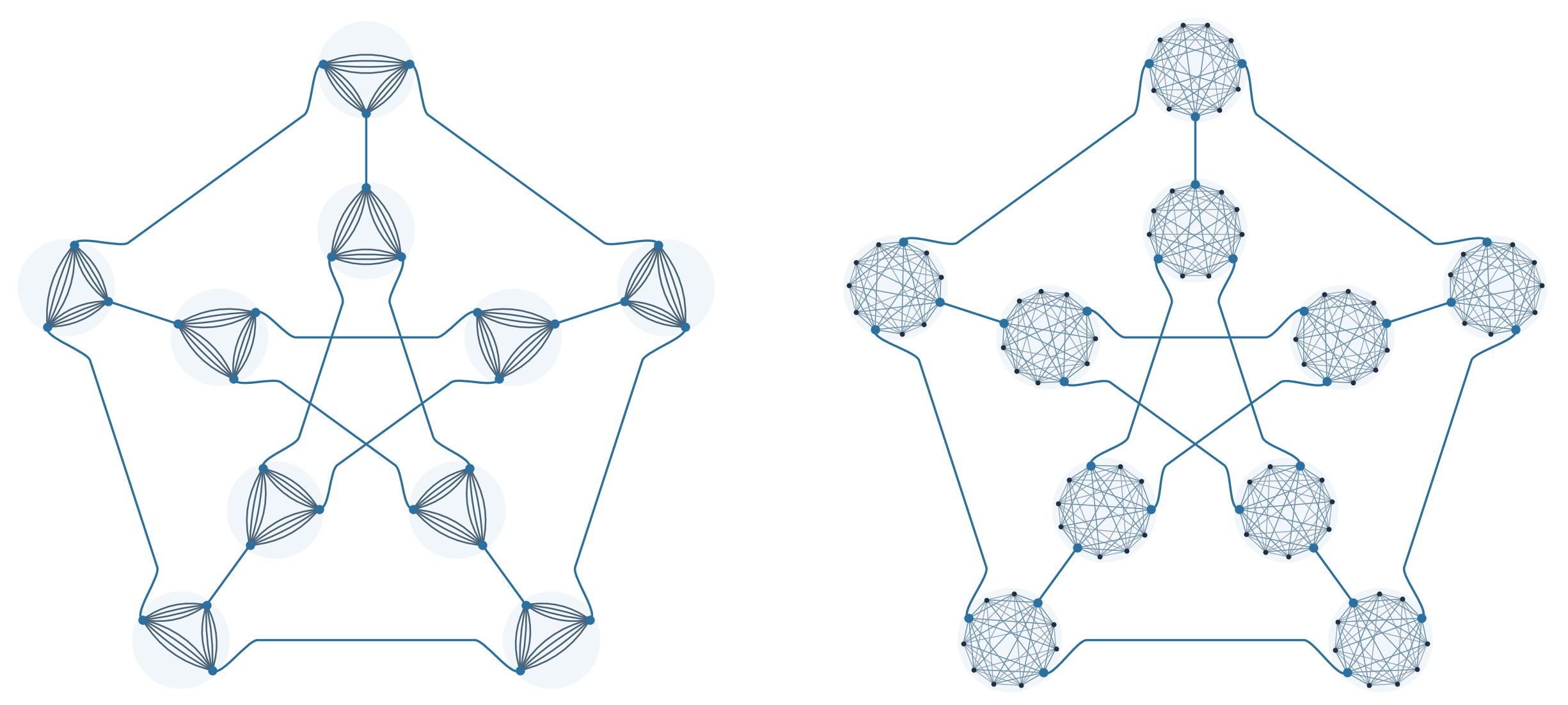}
    \caption{Counterexamples $F_3$ (left) and $G_3$ (right).}
    \label{fig:fig1}
\end{figure}

\begin{figure}[H]
    \centering
    \includegraphics[width=0.66\linewidth]{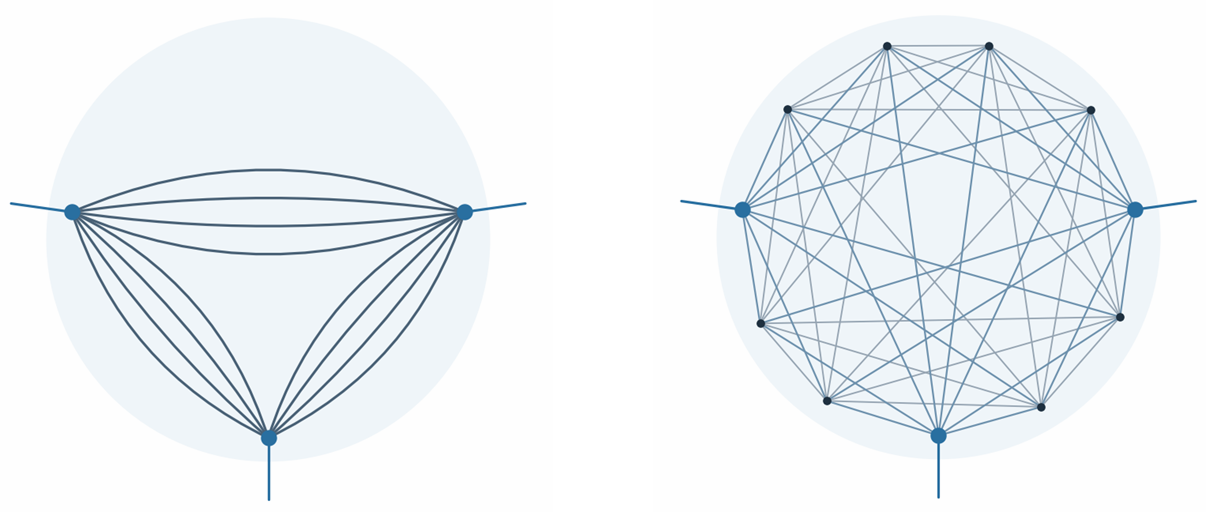}
    \caption{A block of $F_3$ (left) and a block of $G_3$ (right).}
    \label{fig:fig2}
\end{figure}

Figure~\ref{fig:fig1} shows the multigraph $F_3$ and the simple graph $G_3$. Enlarged blocks are shown in Figure~\ref{fig:fig2}. In each block of $F_3$, there are four parallel edges between any two vertices. In each block of $G_3$, the three blue vertices (attachment vertices) are pairwise nonadjacent and joined to all eight black vertices, and the black vertices induce $K_8$ minus a perfect matching.

\section{Proof of Theorem~\ref{thm:connected}}\label{sec:connected}

Let $G$ be a connected graph containing a cycle, and let $\ell\ge2$ be an integer. Fix a spanning tree $T$ of $G$. In any $\ell$-fold cover of $G$, we may relabel the fibers along $T$ so that every tree edge has the identity permutation. In a cyclic cover this relabeling uses translations and makes every tree voltage $0$. Each sheet then contains a copy of $T$, and traversing an edge outside $T$ changes the sheet label by its voltage or its negative. Thus, when all tree voltages are $0$, a cyclic cover is connected if and only if the voltages on the edges outside $T$ generate $\Z_\ell$.

\begin{lemma}\label{lem:connected-min}
Let $G$ be connected and contain a cycle, and let $\ell\ge2$. There is a connected cyclic $\ell$-fold cover of $G$ with chromatic index $s_\ell(G)$, and there is a connected $\ell$-fold cover of $G$ with chromatic index $a_\ell(G)$.
\end{lemma}

\begin{proof}
Choose a cyclic $\ell$-fold cover of minimum chromatic index and relabel its fibers so that its voltages on $T$ are $0$. Let $b$ be the order of the subgroup generated by its voltages. This subgroup consists of the multiples of $\ell/b$ in $\Z_\ell$, and each component corresponds to one of its cosets. Every component is isomorphic to the cyclic $b$-fold cover whose voltages $\beta$ are obtained by dividing the original voltages by $\ell/b$. In particular, the voltages $\beta(e)$ for $e\notin E(T)$ generate $\Z_b$. When $b=1$, each component is a copy of $G$.

Choose voltages $\alpha(e)$ modulo $\ell$ with $\alpha(e)\equiv\beta(e)\pmod b$, taking $\alpha(e)=0$ for $e\in E(T)$. Fix an edge $f\notin E(T)$. For each prime $p$ dividing $\ell$ but not $b$, require $\alpha(f)\equiv1\pmod p$. These choices are possible by the Chinese remainder theorem. For each prime $p$ dividing $b$, at least one $\beta(e)$ is nonzero modulo $p$, and the corresponding $\alpha(e)$ is also nonzero modulo $p$. Thus no prime divisor of $\ell$ divides every voltage $\alpha(e)$, so these voltages generate $\Z_\ell$. Reduction of the sheet labels modulo $b$ gives a projection from the new connected cover onto the cyclic $b$-fold cover with voltages $\beta$. An optimal edge-coloring of that cover lifts to the new cover.

For arbitrary covers, choose an $\ell$-fold cover of minimum chromatic index with an optimal edge-coloring. Fix an edge $f=uv$ on a cycle $C$ of $G$. Each component of the cover projects onto $G$, since any path in $G$ can be lifted from a vertex over its initial vertex. Thus each component contains an edge over $f$. The subgraph consisting of all vertices and edges that project to $C$ is $2$-regular, so every edge over $f$ belongs to a cycle.

If the cover has two distinct components, choose edges $x_1y_1$ and $x_2y_2$ over $f$, one in each component, with $x_1,x_2$ over $u$ and $y_1,y_2$ over $v$. Permute the colors in one component so that both edges have the same color. Replace them by $x_1y_2$ and $x_2y_1$, both projecting onto $f$, and give the new edges the common color. This replacement preserves the matching over $f$ and the proper edge-coloring. Since the deleted edges lie on cycles, the two original components remain connected after the deletions and the new edges join them. Repeating the operation yields a connected cover.
\end{proof}

\begin{lemma}\label{lem:interval}
If $G$ is connected and contains a cycle, then both connected-cover spectra are intervals for every $\ell\ge2$.
\end{lemma}

\begin{proof}
Fix a spanning tree $T$ of $G$ and an edge $f\notin E(T)$. Use identity permutations on the tree edges as above. If the permutation over $f$ is a single cycle of length $\ell$, the cover is connected, regardless of the other permutations. In a cyclic cover we may use translation by $1$ over $f$.

To pass from one connected cover to another of the same type, first replace the permutation over $f$ by this cycle. Then change the other permutations one at a time to their values in the target cover. Finally, change the permutation over $f$ to its target value. Each intermediate cover is connected. At each step the two covers have identical restrictions after the matching over the changed edge is deleted. An edge-coloring of either cover restricts to that common subgraph, and the replacement matching needs at most one new color. Thus consecutive chromatic indices differ by at most $1$, proving that each spectrum is an interval.
\end{proof}

For a graph $J$, an integer $q\ge\Delta(J)$, and a set $W\subseteq V(J)$, define the deficiency
\begin{align}\label{eq:deficiency}
D_q^J(W)
=q|W|-2|E(J[W])|
=|\partial_J(W)|+\sum_{w\in W}\bigl(q-\deg_J(w)\bigr).
\end{align}
An odd set $W$ with at least three vertices has density at most $q$ precisely when $D_q^J(W)\ge q$.

\begin{lemma}\label{lem:critical}
Let $H$ be an edge-critical graph with $\chi'(H)=k\ge\Delta(H)+2$. A cover of $H$ has chromatic index $k$ if and only if one of its components is a $1$-fold cover of $H$.
\end{lemma}

\begin{proof}
Let $q=k-1$. Since $k>\Delta(H)+1$, \eqref{eq:GS} gives an odd vertex set of density greater than $q$. Its induced subgraph has chromatic index at least $q+1$, so criticality implies that the set is $V(H)$. Let $n=|V(H)|$ and $m=|E(H)|$. Then $n$ is odd and $2m>q(n-1)$. Deleting any edge leaves a $q$-edge-colorable graph, so $2(m-1)\le q(n-1)$. Since $q(n-1)$ is even, we have
\begin{align*}
2m=q(n-1)+2,\quad D_q^H(V(H))=q-2.
\end{align*}

Let $J$ be an $\ell$-fold cover with $\chi'(J)=k$. By \eqref{eq:GS} there is an odd set $W\subseteq V(J)$ with $D_q^J(W)<q$. For each $v\in V(H)$, let $x_v$ be the number of vertices of $W$ over $v$, and let $S_i=\{v:x_v\ge i\}$ for $1\le i\le\ell$. As in the counting argument of Chen, van der Holst, Luo, and Ma~\cite[Theorem~2.1]{chen2026spectra}, the matching over a base edge $uv$ contributes at most $\min\{x_u,x_v\}$ edges to $J[W]$. Counting base edges with multiplicity gives
\begin{align*}
|E(J[W])|
\le\sum_{uv\in E(H)}\min\{x_u,x_v\}
=\sum_{i=1}^{\ell}|E(H[S_i])|.
\end{align*}
Consequently,
\begin{align}\label{eq:level}
|W|=\sum_{i=1}^\ell|S_i|,\quad
D_q^J(W)\ge\sum_{i=1}^\ell D_q^H(S_i).
\end{align}
By \eqref{eq:deficiency}, $D_q^H(S_i)\ge|S_i|$ for all $i$, since $q\ge\Delta(H)+1$. For a proper odd set $S_i$, criticality gives a $q$-edge-coloring of $H[S_i]$, so $2|E(H[S_i])|\le q(|S_i|-1)$ and $D_q^H(S_i)\ge q$. This would contradict $D_q^J(W)<q$ by \eqref{eq:level}. Thus every odd $S_i$ equals $V(H)$, and there is an odd number of such sets. Since $q-2=D_q^H(V(H))\ge n$, we have $q\ge3$. Three copies of $V(H)$ would contribute $3(q-2)\ge q$ to the sum in \eqref{eq:level}, so exactly one $S_i$ equals $V(H)$. Any further nonempty $S_i$ has even size and contributes at least $2$, again forcing the sum to be at least $q$.

Since the sets $S_i$ are nested, the only nonempty set is $S_1=V(H)$. Hence $x_v=1$ for every $v\in V(H)$. If even one base edge has no lift in $J[W]$, then $D_q^J(W)\ge qn-2(m-1)=q$, a contradiction. Thus $J[W]$ contains a lift of every base edge, and no edge of $J$ joins $W$ to its complement. The projection maps $J[W]$ isomorphically onto $H$, which is connected by criticality. Conversely, a $1$-fold component forces $\chi'(J)\ge k$, and the coloring of $H$ lifts to give $\chi'(J)\le k$.
\end{proof}

We use the following formula for $s_\ell(H)$, which follows from the critical-graph spectrum theorem of Chen, van der Holst, Luo, and Ma~\cite[Theorem~1.2]{chen2026spectra}. Let $H$ be an edge-critical graph of order $n$ and size $m$ with $\Delta(H)\ge3$ and $\chi'(H)\ge\Delta(H)+1$. For odd $\ell\ge3$,
\begin{align}\label{eq:critical-minimum}
s_\ell(H)=\min\left\{\chi'(H)-1,\max\left\{\Delta(H),\left\lceil\frac{2\ell m}{\ell n-1}\right\rceil\right\}\right\}.
\end{align}

\begin{proof}[Proof of Theorem~\ref{thm:connected}]
Part~(i) follows from Lemmas~\ref{lem:connected-min} and~\ref{lem:interval}.

For part~(ii), assume that $k\ge\Delta(G)+2$. For a nonbridge $e$, assign voltage $1$ to $e$ and voltage $0$ to every other edge. The resulting cyclic cover is connected because $G-e$ contains a spanning tree. It contains a copy of $G-e$ in each sheet. Deleting a single edge reduces the chromatic index by at most $1$, so the upper endpoint of each spectrum is at least $k-1$. If some nonbridge $e$ satisfies $\chi'(G-e)=k$, the corresponding cover attains $k$, the upper bound for every cover of $G$.

Now assume that $\chi'(G-e)=k-1$ for every nonbridge $e$. Choose a subgraph $H$ that is minimal subject to $\chi'(H)=k$. Then $H$ is edge-critical and connected. Every nonbridge $e$ of $G$ belongs to $H$, since otherwise $G-e$ contains $H$. The edges outside $H$ are therefore bridges and form a forest. A path in this forest joining two vertices of $H$, with no internal vertex in $H$, would form a cycle with a path in $H$. Since $G$ is connected, the forest consists of trees attached to $H$ at one vertex each.

Let $J$ be a cover of $G$ with projection $\pi$. Let $J_H$ be the subgraph of $J$ formed by all vertices that project to vertices of $H$ and all edges that project to edges of $H$. Then $\pi$ restricts to a covering projection from $J_H$ onto $H$. The preimage of each attached tree is a disjoint union of copies of that tree, since we can relabel along its edges so that all its permutations are the identity. Thus each lifted tree meets $J_H$ at a single vertex, and $J_H$ is connected whenever $J$ is connected. Moreover,
\begin{align*}
\chi'(J)=\max\{\chi'(J_H),\Delta(G)\}.
\end{align*}
The lower bound follows from $J_H\subseteq J$ and $\Delta(J)=\Delta(G)$. To attain it, extend an edge-coloring of $J_H$ to each tree in order of increasing distance from its attachment vertex. At each step the endpoint farther from the attachment vertex has no incident colored edges, and the other endpoint forbids at most $\Delta(G)-1$ colors. If $J$ is connected, then $J_H$ is a connected $\ell$-fold cover of $H$ with $\ell\ge2$. Since $k\ge\Delta(H)+2$, Lemma~\ref{lem:critical} gives $\chi'(J_H)\le k-1$. The displayed identity then gives $\chi'(J)\le k-1$, proving \eqref{eq:upper-end}.

For part~(iii), let $\ell\ge3$ be odd. If $L_\ell(G)>\Delta(G)$, Theorem~\ref{thm:bound} gives $a_\ell(G)=s_\ell(G)=L_\ell(G)$. In this case $k\ge L_\ell(G)>\Delta(G)$. If $k=\Delta(G)+1$, then $k=L_\ell(G)$, so both spectra are singletons. If $k\ge\Delta(G)+2$, \eqref{eq:upper-end} gives the upper endpoint.

It remains to consider $L_\ell(G)=\Delta(G)$ under the second case of \eqref{eq:upper-end}. Let $H$ be the edge-critical subgraph chosen in that case. Every cyclic cover of $G$ restricts to a cyclic cover of $H$, and every cyclic cover of $H$ extends to one of $G$ by giving the edges of the attached trees voltage $0$. The formula for $\chi'(J)$ therefore gives
\begin{align*}
s_\ell(G)=\max\{\Delta(G),s_\ell(H)\}.
\end{align*}
Since $k\ge\Delta(H)+2$, we have $\Delta(H)\ge3$. Indeed, graphs of maximum degree at most $2$ are disjoint unions of paths and cycles, so their chromatic index is at most their maximum degree plus $1$. Let $q=k-1$, $n=|V(H)|$, and $m=|E(H)|$. The proof of Lemma~\ref{lem:critical} gives $2m=q(n-1)+2$. Since $q\ge3$ and $\ell\ge3$,
\begin{align*}
q(\ell n-1)-2\ell m
&=q(\ell-1)-2\ell\\
&\ge 3(\ell-1)-2\ell=\ell-3\ge0.
\end{align*}
Thus $\left\lceil\frac{2\ell m}{\ell n-1}\right\rceil\le q$. Since $\Delta(H)\le q-1$, formula~\eqref{eq:critical-minimum} gives
\begin{align*}
s_\ell(H)=\max\left\{\Delta(H),\left\lceil\frac{2\ell m}{\ell n-1}\right\rceil\right\}=L_\ell(H).
\end{align*}
For the second equality, $V(H)$ is odd, so the displayed maximum is at most $L_\ell(H)$, which is at most $s_\ell(H)$. Since $H\subseteq G$, we have $L_\ell(G)\ge\max\{\Delta(G),L_\ell(H)\}=s_\ell(G)$. The general lower bound $L_\ell(G)\le a_\ell(G)\le s_\ell(G)$ gives the desired equality.
\end{proof}

For example, let $G$ be a triangle with five parallel edges on each side. Then $\Delta(G)=10$ and $\chi'(G)=15$. The set $S=V(G)$ gives $L_3(G)=\lceil\frac{90}{8}\rceil=12$. Deleting any edge reduces the chromatic index by $1$, so $\mathcal{S}_3^{\mathrm{cyc}}(G)=\mathcal{S}_3^{\mathrm{all}}(G)=[12,14]$. Let $G'$ be obtained by joining a vertex of a $2$-cycle to a vertex of $G$ by a bridge. Then $\Delta(G')=11$. A $15$-edge-coloring of $G$ extends to the three new edges by coloring the bridge first and then the two edges of the $2$-cycle, so $\chi'(G')=15$. Every three-vertex set $S\subseteq V(G')$ with $S\neq V(G)$ induces at most $6$ edges. Since $|E(G')|=18$ and $\Delta(G')=11$, we obtain $L_3(G')=\lceil\frac{90}{8}\rceil=12$, with the maximum attained by $S=V(G)$. Deleting either edge of the $2$-cycle leaves a connected graph of chromatic index $15$. Thus $\mathcal{S}_3^{\mathrm{cyc}}(G')=\mathcal{S}_3^{\mathrm{all}}(G')=[12,15]$.


\section{Polynomial-time algorithms}\label{sec:algorithm}

\begin{proposition}\label{prop:algorithms}
Let $G$ be a graph. For parts~\textup{(ii)} and~\textup{(iii)}, assume that $G$ is connected and contains a cycle.
\begin{enumerate}
\item[\textup{(i)}] For every positive integer $\ell$, the value $L_\ell(G)$ can be computed in polynomial time.
\item[\textup{(ii)}] For every $\ell\ge2$, if $\chi'(G)\ge\Delta(G)+2$, the common upper endpoint in Theorem~\ref{thm:connected}\textup{(ii)} can be computed in polynomial time. In either case of Theorem~\ref{thm:connected}\textup{(iii)}, both endpoints can be computed in polynomial time.
\item[\textup{(iii)}] For fixed odd $\ell\ge3$ with $L_\ell(G)>\Delta(G)$, a connected cyclic $\ell$-fold cover of $G$ with any prescribed chromatic index in the common spectrum can be constructed in polynomial time.
\end{enumerate}
\end{proposition}

\begin{proof}
For part~(i), let $\ell$ be a positive integer. In particular, $L_1(G)=\max\{\Delta(G),\lceil\Gamma(G)\rceil\}$. For an integer $q\ge\Delta(G)$, the inequality $L_\ell(G)\le q$ holds if and only if every odd set $S\subseteq V(G)$ satisfies
\begin{align}\label{eq:threshold}
\ell D_q^G(S)\ge q.
\end{align}
For $|S|\ge3$, \eqref{eq:threshold} is equivalent to $2\ell|E(G[S])|\le q(\ell|S|-1)$, the inequality in the definition of $L_\ell(G)$. A singleton satisfies \eqref{eq:threshold} because its deficiency is $q$.

To test \eqref{eq:threshold}, add a vertex $z$ to $G$. Give each original edge capacity $1$ and join each vertex $v$ of $G$ to $z$ by an edge of capacity $q-\deg_G(v)$. Let $K$ denote the resulting capacitated graph. For $S\subseteq V(G)$, the capacity of $\partial_K(S)$ is $D_q^G(S)$ by \eqref{eq:deficiency}. Let $\mathcal{T}=V(G)$ if $|V(G)|$ is even, and let $\mathcal{T}=V(G)\cup\{z\}$ otherwise. Since $|\mathcal{T}|$ is even, the cuts with an odd number of vertices of $\mathcal{T}$ on either shore are precisely the cuts whose shore avoiding $z$ has odd size. The minimum capacity of such a cut can be found in polynomial time by the odd-cut algorithm of Padberg and Rao~\cite{padberg1982odd}.

Let $T^*$ be a Gomory-Hu tree~\cite{gomory1961flows} of $K$, so $V(T^*)=V(K)$. For an edge $e=uv$ of $T^*$, let $S_e$ be the vertex set of either component of $T^*-e$. The cut $\partial_K(S_e)$ is the fundamental cut associated with $e$. Its capacity equals the weight of $e$ in $T^*$, which is also the minimum capacity of a cut separating $u$ and $v$ in $K$.

Mark $e$ if $|S_e\cap\mathcal{T}|$ is odd. Since $|\mathcal{T}|$ is even, the choice of component defining $S_e$ is immaterial. Let $T'$ be the spanning subgraph of $T^*$ formed by the marked edges. For a vertex $v$ of $T^*$, the components of $T^*-v$ correspond to the edges incident with $v$. Such an edge is marked exactly when the corresponding component contains an odd number of vertices of $\mathcal{T}$. Thus $\deg_{T'}(v)$ has the parity of $|\mathcal{T}\setminus\{v\}|$, which is odd exactly when $v\in\mathcal{T}$.

For a set $S\subseteq V(K)$ with $|S\cap\mathcal{T}|$ odd, the sum of the degrees in $T'$ over $S$ is odd. Edges with both endpoints in $S$ contribute twice, so some marked edge $uv$ has exactly one endpoint in $S$. Its weight is at most the capacity of $\partial_K(S)$, which separates $u$ and $v$. Conversely, the fundamental cut of each marked edge has an odd number of vertices of $\mathcal{T}$ on each shore and has capacity equal to the weight of that edge. Therefore the minimum weight of a marked edge equals the minimum capacity of such a cut.

For a graph with at least one edge, binary search over the integers from $\Delta(G)$ to $|E(G)|$ computes $L_\ell(G)$. An edgeless graph has $L_\ell(G)=0$.

For part~(ii), by \eqref{eq:GS}, the condition $\chi'(G)\ge\Delta(G)+2$ is equivalent to $L_1(G)\ge\Delta(G)+2$, and in that case $k=\chi'(G)=L_1(G)$. Also, $\chi'(G-e)=k$ if and only if $L_1(G-e)=k$. One implication uses $L_1(G-e)\le\chi'(G-e)\le k$, and the other uses \eqref{eq:GS} and $\Delta(G-e)+1<k$. Thus computing $L_1(G-e)$ for each nonbridge determines the upper endpoint in \eqref{eq:upper-end}.

Now let $\ell\ge3$ be odd. If $L_\ell(G)>\Delta(G)$, then $L_1(G)\ge L_\ell(G)>\Delta(G)$, so $k=L_1(G)$. The common lower endpoint is $L_\ell(G)$. The preceding procedure gives the upper endpoint when $k\ge\Delta(G)+2$, and the spectrum is $\{k\}$ when $k=\Delta(G)+1$. If $L_\ell(G)=\Delta(G)$, $k\ge\Delta(G)+2$, and all nonbridge deletions give $\chi'(G-e)=k-1$, the common endpoints are $L_\ell(G)$ and $k-1$ by Theorem~\ref{thm:connected}\textup{(iii)}. The tests of $L_1(G-e)$ also identify this case.

For part~(iii), Lemma~\ref{lem:pairs} gives a cyclic cover attaining the lower endpoint. The cyclic construction in Lemma~\ref{lem:connected-min} makes it connected using only its voltage assignment and the subgroup generated by the voltages after relabeling along a spanning tree. The coloring argument in that lemma certifies the bound on the chromatic index of the new cover. For fixed $\ell$, all these voltage operations take polynomial time.

When the spectrum is a singleton, the first cover suffices. Otherwise, use voltage $1$ on a suitable nonbridge and voltage $0$ elsewhere to attain the upper endpoint, as in the proof of Theorem~\ref{thm:connected}. Lemma~\ref{lem:interval} gives a sequence between the two covers with at most $|E(G)|+2$ steps. For every cover $J$ in this sequence, taking full preimages of odd vertex sets of $G$ gives $L_1(J)\ge L_\ell(G)>\Delta(J)$. By \eqref{eq:GS}, $\chi'(J)=L_1(J)$, so computing $L_1(J)$ along the sequence finds a cover with the prescribed index. The procedure outputs the cover and certifies its chromatic index.
\end{proof}

\section{A concluding question}
For odd $\ell\ge3$, Theorem~\ref{thm:bound} shows that $s_\ell(G)\neq a_\ell(G)$ is possible only if $a_\ell(G)=L_\ell(G)=\Delta(G)$ and $s_\ell(G)=\Delta(G)+1$. For the graphs in Corollary~\ref{cor:counterexamples}, we have $s_\ell(G)=a_\ell(G)=\Delta(G)+1$. Does $s_\ell(G)=a_\ell(G)$ hold for every loopless graph $G$ and every odd integer $\ell\ge3$?

\section*{Declaration of generative AI use}
During the preparation of this work, the author used OpenAI's GPT-6.1 Sol to find the Petersen counterexamples in Corollary~\ref{cor:counterexamples}. The author independently reviewed, verified, and edited all AI-generated content and takes full responsibility for the final manuscript.


\begin{thebibliography}{10}

\bibitem{chen2024short}
G.~Chen, Y.~Hao, X.~Yu, and W.~Zang.
\newblock A short proof of the {G}oldberg-{S}eymour conjecture.
\newblock {\em arXiv preprint arXiv:2407.09403}, 2024.

\bibitem{chen2025goldberg}
G.~Chen, G.~Jing, and W.~Zang.
\newblock {P}roof of the {G}oldberg-{S}eymour conjecture on edge-colorings of multigraphs.
\newblock {\em J. Comb. Optim.}, 50(3):Paper No. 23, 2025.

\bibitem{chen2026spectra}
G.~Chen, H.~van~der Holst, R.~Luo, and Y.~Ma.
\newblock On the spectra of chromatic number and chromatic index of cyclic covers.
\newblock {\em arXiv preprint arXiv:2608.02934}, 2026.

\bibitem{gomory1961flows}
R.~E. Gomory and T.~C. Hu.
\newblock Multi-terminal network flows.
\newblock {\em J. Soc. Ind. Appl. Math.}, 9(4):551--570, 1961.

\bibitem{gross1974voltage}
J.~L. Gross.
\newblock Voltage graphs.
\newblock {\em Discrete Math.}, 9:239--246, 1974.

\bibitem{gross1977coverings}
J.~L. Gross and T.~W. Tucker.
\newblock Generating all graph coverings by permutation voltage assignments.
\newblock {\em Discrete Math.}, 18(3):273--283, 1977.

\bibitem{luo2026flows}
R.~Luo, Y.~Ma, D.~C. Stephens, and D.~Ye.
\newblock Integer flows on cyclic covers of graphs.
\newblock {\em arXiv preprint arXiv:2610.29125}, 2026.

\bibitem{padberg1982odd}
M.~W. Padberg and M.~R. Rao.
\newblock Odd minimum cut-sets and {$b$}-matchings.
\newblock {\em Math. Oper. Res.}, 7(1):67--80, 1982.

\bibitem{pisanski1983bundles}
T.~Pisanski, J.~Shawe-Taylor, and J.~Vrabec.
\newblock Edge-colorability of graph bundles.
\newblock {\em J. Combin. Theory Ser. B}, 35(1):12--19, 1983.

\bibitem{plachta2021coverings}
L.~Plachta.
\newblock Coverings of cubic graphs and 3-edge colorability.
\newblock {\em Discuss. Math. Graph Theory}, 41(1):311--334, 2021.

\end{thebibliography}
\end{document}